\documentclass[english,10pt]{article}
\usepackage{amsmath,amssymb,amsbsy,graphicx, setspace,cite}
\usepackage{color,lpic}
\title{Suppression of Limit Cycle Oscillations \\ using the Nonlinear Tuned Vibration Absorber}

\begin{document}
\date{}
\maketitle

\begin{center}
{This is an author post-print, original manuscript at\\
http://rspa.royalsocietypublishing.org/content/471/2176/20140976
\\
G. Habib, G. Kerschen\\\vspace{0.8cm}

\small Space Structures and Systems Laboratory\\
Department of Aerospace and Mechanical Engineering\\
University of Li\`ege, Li\`ege, Belgium \\
E-mail: giuseppe.habib,g.kerschen@ulg.ac.be\\\vspace{0.5cm} \vspace{1cm}

\vspace{-0.5cm}
\begin{abstract}
The objective of the present study is to mitigate, or even completely eliminate, the limit cycle oscillations in mechanical systems using a passive nonlinear absorber, termed the nonlinear tuned vibration absorber (NLTVA). An unconventional aspect of the NLTVA is that the mathematical form of its restoring force is not imposed a priori, as it is the case for most existing nonlinear absorbers. The NLTVA parameters are determined analytically using stability and bifurcation analyses, and the resulting design is validated using numerical continuation. The proposed developments are illustrated using a Van der Pol-Duffing primary system.

\vspace{1cm}

\noindent \emph{Keywords}: limit cycle oscillations, vibration absorber, nonlinearity synthesis, stability analysis, bifurcation analysis.
\end{abstract}

\vspace{1.5cm} Corresponding author: \\ Giuseppe Habib\\
Space Structures and Systems Laboratory\\
Department of Aerospace and Mechanical Engineering\\
University of Li\`ege
\\ 1 Chemin des Chevreuils (B52/3), B-4000 Li\`ege, Belgium. \\
Email: giuseppe.habib@ulg.ac.be.
\vspace{2cm}\\}
\end{center}

\section{Introduction}

Limit cycle oscillations (LCOs) are encountered in a number of real engineering applications including aircraft \cite{Denegri,Trickey}, machine tools \cite{Stepan,Mann} and automotive disk brakes \cite{Hagedorn}. LCOs often limit the performance and can also endanger the safety of operation \cite{Griffin}.

Active control strategies have been proposed as a means of counteracting LCOs \cite{Strganac,Cooper,vandeWouw}. These references have shown that active control can be used to raise the threshold above which LCOs occur. However, active control is also limited by its requirements in terms of energy or space for actuators. Furthermore, delay in the feedback loop can generate unexpected instabilities \cite{Stepan2} whereas saturation of the actuators can limit the robustness of stability \cite{Habib}.

Passive vibration absorbers represent another alternative for mitigating undesired LCOs. Specifically, the linear tuned vibration absorber (LTVA), which comprises a small mass attached to the host system through a damper and a spring, has been widely studied in the literature \cite{Mansour,Rowbottom, Fujino, Gattulli, Gattulli2}. In most of these works, the system under investigation is the classical Van der Pol (VdP) oscillator. References \cite{Rowbottom,Fujino} provide simple rules to properly tune the LTVA parameters, while references \cite{Gattulli,Gattulli2} study the post-bifurcation behavior of the coupled system.
In \cite{Gattulli3}, a nonlinear damping element was added in parallel with the LTVA to decrease the maximum LCO amplitude. Other nonlinear vibration absorbers, including the autoparametric vibration absorber \cite{Haxton,Vyas}, the nonlinear energy sink (NES) \cite{Lee,Gendelman,Luongo,Lee1,Lee2,Lee3,Lee4,Gendelman2,Tumur} and the hysteretic tuned vibration absorber \cite{Lacarbonara}, have also been considered to increase the effectiveness of vibration attenuation.
In particular, the NES exhibited three different mechanisms to suppress LCOs, namely complete, partial and intermittent suppression. These mechanisms were studied numerically \cite{Lee1}, experimentally \cite{Lee2} and analytically \cite{Gendelman2}.

The main idea of this study is to utilize the nonlinear tuned vibration absorber (NLTVA) for LCO suppression. This absorber, first introduced in \cite{Habib2}, possesses a linear spring and a nonlinear spring whose mathematical form is determined according to the nonlinearity in the host system.
Following references \cite{Mansour,Rowbottom, Fujino, Gattulli, Gattulli2}, the linear spring coefficient is determined to maximize the stable region of the trivial solution of the host system. Subsequently, the nonlinear spring is designed to ensure supercritical behavior and to mitigate the LCOs in the postcritical range. A fundamental result of this paper is that, if properly designed, the linear and nonlinear springs of the NLTVA can complement each other giving rise to a very effective LCO suppression and management strategy.
The example that will serve to validate the proposed developments is the Van der Pol--Duffing (VdPD) oscillator \cite{Guckenheimer}, which is a paradigmatic model for the description of self-excited oscillations.
The bifurcation behavior of the VdPD oscillator was studied in \cite{Szemplinska} whereas its stabilization using active control was proposed in \cite{Xu,Li}.

The paper is organized as follows. Section \ref{ProblemStatement} introduces the three design objectives pursued in this paper. In Section \ref{Elimination}, optimal values for the linear parameters of the NLTVA are determined using stability analysis of the coupled system. Section \ref{BifurcationAnalysis} investigates the bifurcations occurring at the loss of stability and proposes an analytical tuning rule for the nonlinear coefficient of the NLTVA. In Section \ref{Reduction}, the reduction of the LCO amplitude in the postcritical range is discussed.
In Section \ref{Validation}, the analytical results are validated numerically using the MATCONT software and the global behavior of the system is also carefully discussed.
The NLTVA is then compared with the NES, highlighting the better performance of the former absorber. Finally, conclusions are drawn in Section \ref{Conclusions}.

\section{Problem Formulation}\label{ProblemStatement}

The primary system considered throughout this work is the Van der Pol--Duffing (VdPD) oscillator:
\begin{equation}
m_1 q''_1+c_1\left(q_1^2-1\right)q'_1+k_1q_1+k_{nl1}q_1^3=0
\end{equation}
where $m_1$, $c_1$, $k_1$ and $k_{nl1}$ are the oscillator's mass, damping and the coefficients of the linear and cubic springs, respectively. For instance, for an in-flow wing, the terms $c_1\left(q_1^2-1\right)q'_1$ and $k_{nl1}q_1^3$ would model the fluid-structure interaction and the structural nonlinearity, respectively.
Stability analysis demonstrates that the trivial equilibrium point of the system loses stability when $\mu_1=c_1/(2\sqrt{k_1m_1})=0$. Loss of stability occurs through either a supercritical or a subcritical Hopf bifurcation. This latter scenario is dangerous, because stable large-amplitude LCOs can co-exist with the stable equilibrium point \cite{NayfehBalachandran}.

The objective of the present study is to mitigate, or even completely eliminate, the LCOs of the VdPD oscillator through the attachment of a fully passive nonlinear vibration absorber, termed the NLTVA \cite{Habib2}. One salient feature of the NLTVA compared to existing nonlinear absorbers is that the absorber's load-deflection curve is not imposed a priori, but it is rather synthesized according to the nonlinear restoring force of the primary system. The equations of motion of the coupled VdPD and NLTVA system are:
\begin{eqnarray}
&&m_1 q''_1+c_1\left(q_1^2-1\right)q'_1+k_1q_1+ k_{nl1}q_1^3+c_2\left(q'_1-q'_2\right)+G\left(q_d\right)=0\nonumber\\
&&m_2q''_2+c_2\left(q'_2-q'_1\right)-G\left(q_d\right)=0
\end{eqnarray}
where $m_2$ and $c_2$ are the absorber's mass and viscous damping, respectively. The NLTVA is assumed to have a generic smooth elastic force $G(q_d)$ where $q_d=q_1-q_2$ with $G(0)=0$.

\begin{figure}
\begin{centering}
\includegraphics[trim = 10mm 10mm 10mm 10mm,width=0.32\textwidth]{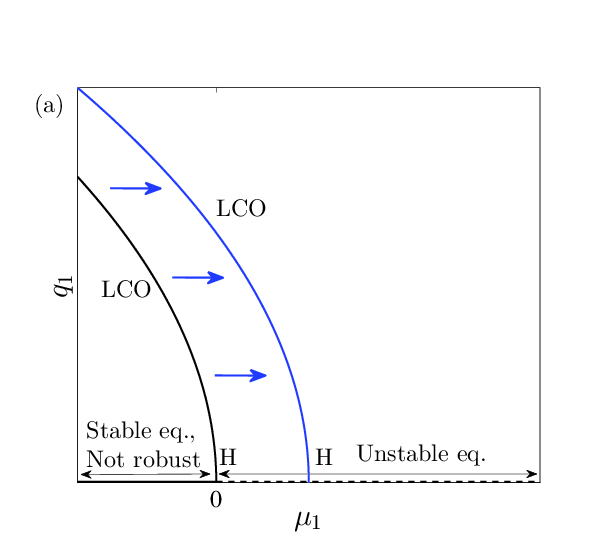}
\includegraphics[trim = 10mm 10mm 10mm 10mm,width=0.32\textwidth]{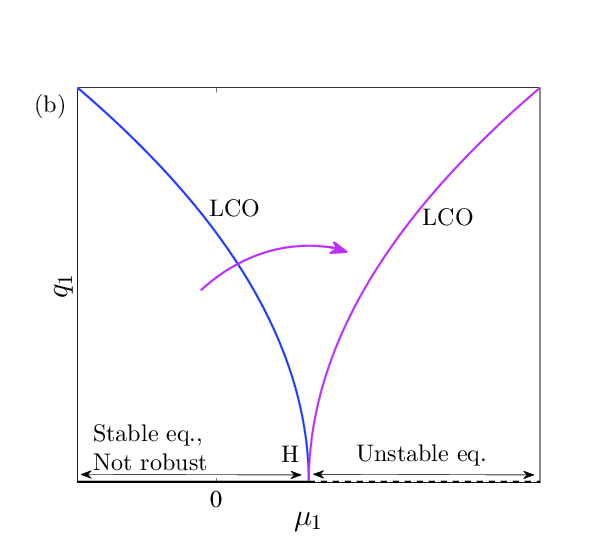}
\includegraphics[trim = 10mm 10mm 10mm 10mm,width=0.32\textwidth]{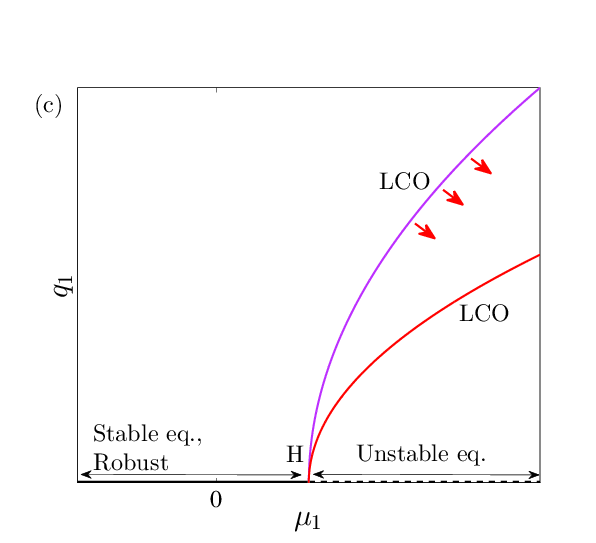}
\par\end{centering}
\caption{\label{explanation}Subcritical VdPD oscillator with an attached NLTVA. (a) Enlargement of the region where the equilibrium point of the VdPD oscillator is stable; (b) transformation of the subcritical Hopf bifurcation into a supercritical Hopf bifurcation; (c) reduction of the amplitude of the remaining LCOs. H and LCO stands for Hopf bifurcation and limit cycle oscillation, respectively. The other labels (stable, unstable, robust) refer to the situation before the introduction of the NLTVA.}
\end{figure}

The design problem is as follows. The mass ratio $\varepsilon=m_2/m_1$ (and, hence, the absorber mass) is prescribed by obvious practical constraints; $\varepsilon=0.05$ is considered in the numerical examples of this paper. The damping coefficient $c_2$ and the absorber's stiffness $G(q_d)$ should be determined so as to:
\begin{enumerate}
  \item Maximize the region where the trivial equilibrium is stable by displacing the Hopf bifurcation toward large positive values of $\mu_1$ (Fig. \ref{explanation}(a));
  \item Avoid a catastrophic bifurcation scenario by transforming the potentially subcritical Hopf bifurcation of the VdPD into a supercritical Hopf bifurcation of the coupled system (Fig. \ref{explanation}(b));
  \item Reduce the amplitude of the remaining LCOs (Fig. \ref{explanation}(c)).
\end{enumerate}
These three design objectives are studied in detail in the next three sections.

\section{Maximization of the stable region of the trivial solution}\label{Elimination}

The first design objective is to stabilize the trivial solution of the VdPD oscillator for values of $\mu_1$ as great as possible. Because the stability of an equilibrium point of a nonlinear system is governed only by the local underlying linear system, the NLTVA should comprise a linear spring for increased flexibility, i.e., $G\left(q_d\right)=k_2\left(q_d\right)+G_{nl}\left(q_d\right)$. The system of interest for the stability analysis is therefore
\begin{eqnarray}
&&m_1 q''_1-c_1q'_1+k_1q_1+c_2\left(q'_1-q'_2\right)+k_2\left(q_1-q_2\right)=0\nonumber\\
&&m_2q''_2+c_2\left(q'_2-q'_1\right)+k_2\left(q_2-q_1\right)=0.
\end{eqnarray}

Introducing the parameters $\omega_{n1}^2=k_1/m_1$, $\omega_{n2}^2=k_2/m_2$, $\mu_2=c_2/(2m_2\omega_{n2})$, $\gamma=\omega_{n2}/\omega_{n1}$ and the dimensionless time $\tau=t\omega_{n1}$, we transform the system into first-order differential equations
\begin{equation}
\left [\begin{array}{c}
\dot x_1\\\dot x_2\\\dot x_3\\\dot x_4
\end{array}\right ]=\left [\begin{array}{cccc}
0&1&0&0\\
-1&2\mu_1&-\gamma^2\varepsilon &-2\mu_2\gamma\varepsilon\\
0&0&0&1\\
-1&2\mu_1&-\gamma^2\left(1+\varepsilon\right)&-2\mu_2\gamma\left(1+\varepsilon\right)
\end{array}\right]\left[\begin{array}{c}
x_1\\x_2\\x_3\\x_4
\end{array}\right]\label{in_matrix}
\end{equation}
or in compact form $\dot{\mathbf x}=\mathbf{Wx}$, where $x_1=q_1$, $x_2=\dot q_1$, $x_3=q_d$ and $x_4=\dot q_d$. The dot indicates derivation with respect to the dimensionless time $\tau$.

\begin{figure}[t]
\begin{centering}
\includegraphics[trim = 10mm 10mm 10mm 10mm,width=0.4\textwidth]{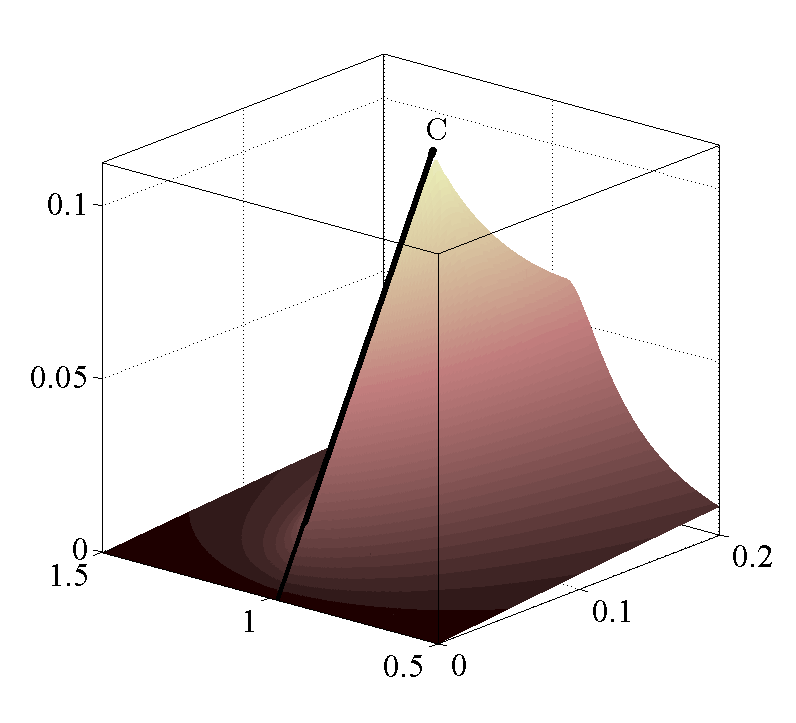}
\put(-210,80){$\mu_1$}
\put(-35,3){$\mu_2$}
\put(-150,5){$\gamma$}
\par\end{centering}
\caption{\label{stab_3d}Stability chart in the $\mu_1,\mu_2,\gamma$ space for $\varepsilon=0.05$. The surface indicates the stability boundary.}
\end{figure}

As reported for a similar system in \cite{Gattulli}, the trivial solution of Eq. (\ref{in_matrix}) is asymptotically stable if and only if the roots of the characteristic polynomial $\det\left(\mathbf W-z\mathbf I\right)=0$ have negative real parts. The roots are computed by solving
$z^4+2 \left((\varepsilon+1) \gamma \mu_2-\mu_1\right)z^3+ \left((\varepsilon+1) \gamma^2-4 \gamma \mu_1 \mu_2+1\right)z^2+2 \gamma (\mu_2-\gamma \mu_1)z+\gamma^2=0$,
which is rewritten as $a_4z^4+a_3z^3+a_2z^2+a_1z+a_0=0$.
Fig. \ref{stab_3d} depicts the stability chart in the $\mu_1,\mu_2,\gamma$ space obtained from direct evaluation of the roots. The surface, which represents the stability boundary, peaks along the $\mu_1$ axis at point C, meaning that the trivial solution can no longer be stable beyond this point. The corresponding maximum value of $\mu_1$ is around 0.1, and the other parameters are $\mu_2 \approx 0.1$ and $\gamma \approx 1$.

In order to calculate analytically the coordinates of point C, thus defining optimal parameters for maximizing stability, the Routh-Hurwitz stability criterion is exploited. The characteristic polynomial has roots with negative real parts if and only if the coefficients $a_i>0,\,i=1,...,4$, $e_2=\left(a_3a_2-a_4a_1\right)/a_3>0$ and $e_3=c_2a_1-a_3a_0>0$, i.e.,
\begin{eqnarray}
a_4&=&1>0\label{a4}\\
a_3&=&2\left((\varepsilon+1) \gamma \mu_2-\mu_1\right)>0\,\iff\,\mu_1<(\varepsilon+1) \gamma \mu_2\\
a_2&=&(\varepsilon+1) \gamma^2-4 \gamma \mu_1 \mu_2+1>0\,\iff\,\mu_1<\frac{(\varepsilon+1) \gamma^2+1}{4 \gamma\mu_2}\\
a_1&=& 2 \gamma (\mu_2-\gamma \mu_1)>0\,\iff\,\mu_1<\frac{\mu_2}{\gamma}\\
a_0&=&\gamma^2>0\\
e_2&=&\frac{2 \left((\varepsilon+1) \gamma \mu_2-\mu_1\right) \left((\varepsilon+1) \gamma^2-4 \gamma \mu_1 \mu_2+1\right)+2 \gamma (\gamma \mu_1-\mu_2)}{2 ((\varepsilon+1) \gamma \mu_2-\mu_1)}>0\\
e_3&=&\frac{\gamma (\mu_2-\gamma \mu_1) \left(2 \left((\varepsilon+1) \gamma \mu_2-\mu_1\right) \left((\varepsilon+1) \gamma^2-4 \gamma \mu_1 \mu_2+1\right)+2 \gamma (\gamma \mu_1-\mu_2)\right)}{(\varepsilon+1) \gamma \mu_2-\mu_1}\nonumber\\
&&\quad-2 \gamma^2 ((\varepsilon+1) \gamma \mu_2-\mu_1)\label{c3}>0.
\end{eqnarray}

\begin{figure}[!]
\begin{centering}
\includegraphics[trim = 10mm 10mm 8mm 10mm,width=0.32\textwidth]{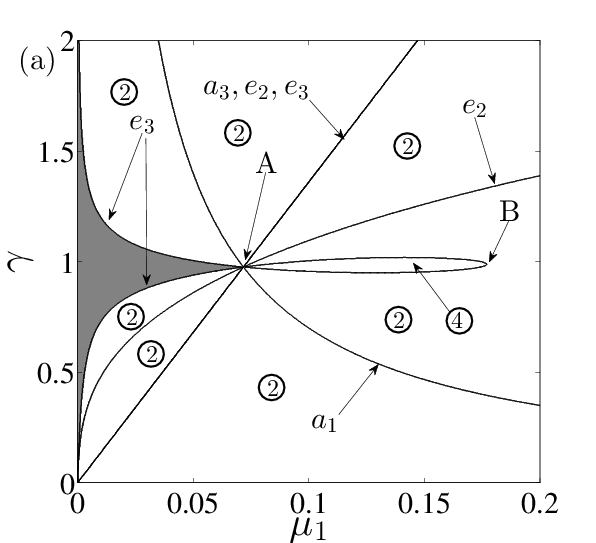}
\includegraphics[trim = 10mm 10mm 8mm 10mm,width=0.32\textwidth]{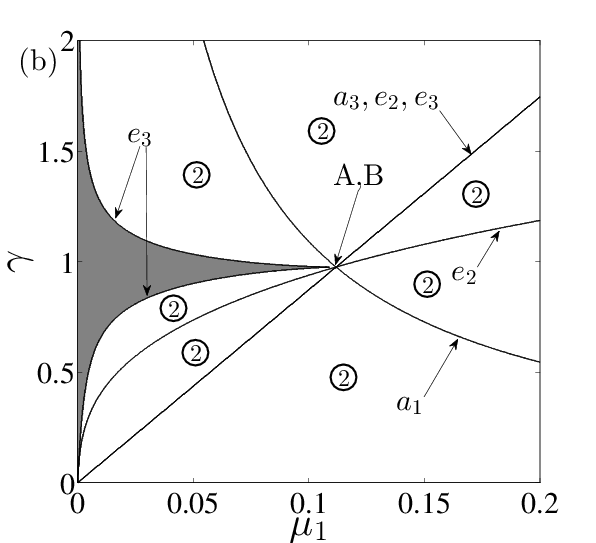}
\includegraphics[trim = 10mm 10mm 8mm 10mm,width=0.32\textwidth]{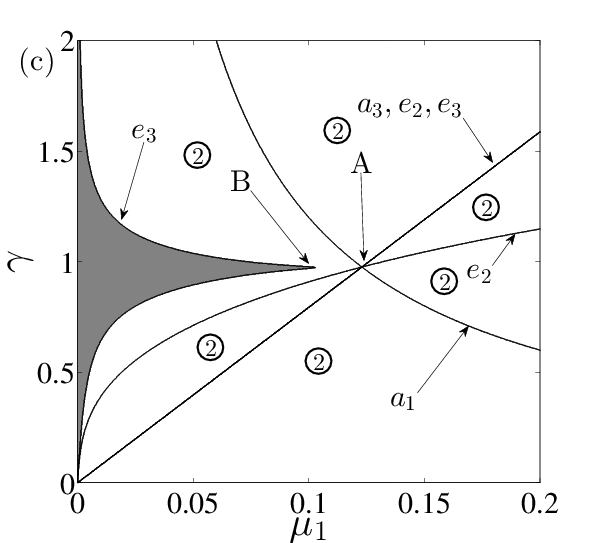}
\par\end{centering}
\caption{\label{mu1_lambda}Sections of the stability chart for different values of $\mu_2$. (a) $\mu_2=0.07$; (b) $\mu_2=1/2\sqrt{\varepsilon/(1+\varepsilon)}=0.1091$; (c) $\mu_2=0.12$. Shaded area: stable, clear area: unstable. Encircled numbers 2 and 4 indicate regions with 2 or 4 eigenvalues (two pairs of complex conjugate eigenvalues) with positive real parts.}
\end{figure}

It is not trivial to interpret the physical meaning of these coefficients, however, they give important information regarding stability.
Fig. 3 represents 2D sections of the 3D stability chart for different values of $\mu_2$. The curves $a_3=0$ and $a_1=0$ intersect at point A, whose coordinates in the ($\mu_1,\gamma$) space are $\text A=\left(\mu_2\sqrt{1+\varepsilon},1/\sqrt{1+\varepsilon}\right)$. Substituting $\gamma=1/\sqrt{1+\varepsilon}$ in $e_3=0$, point $\text B=\left(\varepsilon/\left(4\mu_2\sqrt{1+\varepsilon}\right),1/\sqrt{1+\varepsilon}\right)$ is defined.
As presented in Figs. \ref{mu1_lambda}(a) and (c), A and B mark alternatively the maximal value of $\mu_1$ for stability. Points A and B coincide if and only if $\mu_2=1/2\sqrt{\varepsilon/(1+\varepsilon)}$, which is therefore the optimal condition for stability; this scenario is illustrated in Fig. \ref{mu1_lambda}(b).

Summarizing, optimal values of the linear parameters are
\begin{equation}
\gamma_{\text{opt}}=\frac{1}{\sqrt{1+\varepsilon}},\quad \mu_{2\text{opt}}=\frac{1}{2}\sqrt{\frac{\varepsilon}{1+\varepsilon}}\label{lambdaopt}
\end{equation}
and the corresponding maximal value of $\mu_1$ that guarantees stability is
\begin{equation}
\mu_{1\text{max}}=\frac{\sqrt{\varepsilon}}{2}.\label{max_mu1}
\end{equation}

\section{Enforcement of Supercritical Hopf Bifurcations through Normal Form Analysis}\label{BifurcationAnalysis}

The second design objective is to ensure the robustness of the trivial solution, i.e., no stable LCO can coexist with the stable equilibrium, as depicted in Fig. \ref{explanation}(b). Since supercritical Hopf bifurcations are sought in the coupled system, a detailed investigation of the bifurcations occurring at the loss of stability is the main focus of the present section. Because bifurcation characterization depends on the nonlinear coefficient of the NLTVA, this analysis will allow us to define the optimal value of this coefficient whereas the linear coefficients of the NLTVA should remain close to their optimal values (\ref{lambdaopt}), i.e., $\gamma=0.976$ and $\mu_2=0.109$ for $\varepsilon=0.05$.

Another key element that remains to be determined is the mathematical expression of the NLTVA's elastic force $G(q_d)$.
Applying the 'principle of similarity' between the primary and secondary systems, first proposed for linear systems \cite{dellIsola,Luongo2} and extended to forced nonlinear vibrations in \cite{Habib2}, the mathematical form of the nonlinear spring of the NLTVA is chosen to be also cubic.
Hence, the coupled system writes
\begin{eqnarray}
&&m_1 q''_1+c_1\left(q_1^2-1\right)q'_1+k_1q_1+ k_{nl1}q_1^3+c_2\left(q'_1-q'_2\right)+k_2\left(q_1-q_2\right)+k_{nl2}\left(q_1-q_2\right)^3=0\nonumber\\
&&m_2q''_2+c_2\left(q'_2-q'_1\right)+k_2\left(q_2-q_1\right)+k_{nl2}\left(q_2-q_1\right)^3=0.
\end{eqnarray}
Considering dimensionless coordinates, we transform the system into first-order differential equations
\begin{equation}
\left [\begin{array}{c}
\dot x_1\\\dot x_2\\\dot x_3\\\dot x_4
\end{array}\right ]=\left [\begin{array}{cccc}
0&1&0&0\\
-1&2\mu_1&-\gamma^2\varepsilon &-2\mu_2\gamma\varepsilon\\
0&0&0&1\\
-1&2\mu_1&-\gamma^2\left(1+\varepsilon\right)&-2\mu_2\gamma\left(1+\varepsilon\right)
\end{array}\right]\left[\begin{array}{c}
x_1\\x_2\\x_3\\x_4
\end{array}\right]+\left[\begin{array}{c}
0\\-2\mu_1x_1^2x_2-\alpha_3x_1^3-\beta_3\varepsilon x_3^3\\0\\-2\mu_1x_1^2x_2-\alpha_3x_1^3-\beta_3(1+\varepsilon)x_3^3
\end{array}\right] \label{inx}
\end{equation}
or in compact form $\dot{\mathbf{x}}=\mathbf{Wx}+\mathbf b$. The variables $\alpha_3=k_{nl1}/k_1$ and $\beta_3=k_{nl2}/(k_1\varepsilon)$ have been introduced in these equations.

When stability is lost, one or two pairs of complex conjugate eigenvalues of $\mathbf{W}$ leave the left half plane, which corresponds to single or double Hopf bifurcation, respectively.
As reported in \cite{Gattulli}, a double Hopf bifurcation is likely to occur when $\gamma = \gamma_{\text{opt}}$ and $\mu_2\leq \mu_{\text{opt}}$.
This assertion is confirmed in Fig. \ref{mu1_lambda}(a), where it is shown that losing stability through point A, the system winds up in a region with two pairs of complex conjugate eigenvalues with positive real part.
Referring to the 3D stability chart, the locus of double Hopf bifurcations is depicted by the black line in Fig. \ref{stab_3d}.
The analysis of the double Hopf bifurcation is of practical importance because it can generate quasiperiodic solutions, which might compromise the robustness of the stable trivial solution.
However, a detailed study of the double Hopf bifurcations is beyond the scope of this paper.

\subsection{Single Hopf bifurcation}\label{SHBIF}

The analysis is first focused on the single Hopf bifurcation scenario for which $\mathbf{W}$ has a pair of complex conjugate eigenvalues with zero real part $\lambda_{1,2}=k_1\pm j\omega_1$ and two other eigenvalues $\lambda_3$ and $\lambda_4$ with negative real parts. Vectors $\mathbf s_1$, $\mathbf s_2$, $\mathbf s_3$ and $\mathbf s_4$ are the corresponding eigenvectors. In order to decouple the linear part of the system, we define the transformation matrix
\begin{equation}
\mathbf T=\Bigg[\begin{array}{cccc}
\text{Re}(\mathbf s_1)&\text{Im}(\mathbf s_1)&\text{Re}(\mathbf s_3)&\text{Im}(\mathbf s_3)
\end{array}\Bigg]\label{T}
\end{equation}
where, if $\lambda_3$ and $\lambda_4$ are real, $\text{Re}(\mathbf s_3)$ and $\text{Im}(\mathbf s_3)$ are substituted with $\mathbf s_3$ and $\mathbf s_4$, respectively. Applying the transformation $\mathbf x=\mathbf{Ty}$, Eq. (\ref{inx}) becomes
\begin{equation}
\dot{\mathbf y}=\mathbf{Ay}+\mathbf T^{-1}\mathbf b\label{iny1}
\end{equation}
where \begin{equation}
\mathbf A=\mathbf T^{-1}\mathbf{WT}=\left[\begin{array}{cccc}k_1&\omega_1&0&0\\-\omega_1&k_1&0&0\\0&0&\text{Re}(\lambda_3)&\text{Im}(\lambda_3)\\0&0&\text{Im}(\lambda_4)&\text{Re}(\lambda_4)\end{array}\right].
\end{equation}
The linear part of Eq. (\ref{iny1}) is therefore decoupled.

Variables related to the bifurcation are separated from the other variables using center manifold reduction \cite{NayfehBalachandran}. Relations $y_3=h_3(y_1,y_2)$ and $y_4=h_4(y_1,y_2)$ are used to reduce the dimension of the system to two, while keeping the local dynamics intact. However, since the system has only cubic nonlinearities, $h_3$ and $h_4$ have only cubic and higher-order terms. Thus, for our purpose, they can be considered identically equal to zero, i.e., $y_3\approx0$ and $y_4\approx0$. System (\ref{iny1}) reduces to
\begin{equation}
\left [\begin{array}{c}
\dot y_1\\\dot y_2
\end{array}\right ]=\left[\begin{array}{cc}k_1&\omega_1\\-\omega_1&k_1\end{array}\right]\left[\begin{array}{c}
y_1\\y_2 \end{array}\right]+\left[\begin{array}{c}
d_{130}y_1^3+d_{121}y_1^2y_2+d_{112}y_1y_2^2+d_{103}y_2^3\\
d_{230}y_1^3+d_{221}y_1^2y_2+d_{212}y_1y_2^2+d_{203}y_2^3
\end{array}\right].\label{after_CMR}
\end{equation}

Performing several transformations, namely transformation in complex form, near-identity transformation and transformation in polar coordinates, the bifurcation can be characterized through its normal form
\begin{equation}
\dot r=k_1r+\delta r^3 \label{singH_nf}
\end{equation}
where $\delta=\left(3d_{130}+d_{112}+d_{221}+3d_{203}\right)/8$.
Details of this standard procedure can be found in \cite{Kuznetsov}.
Eq. (\ref{singH_nf}) has solutions $r_0=0$ and $r^*=\sqrt{-k_1/\delta}$. The coefficient $\delta$ can be expressed as a linear function of the nonlinear coefficients $\alpha_3$ and $\beta_3$:
\begin{equation}
\delta=\delta_0(\varepsilon,\gamma,\mu_1,\mu_2)+\delta_{\alpha}(\varepsilon,\gamma,\mu_1,\mu_2)\alpha_3+\delta_{\beta}(\varepsilon,\gamma,\mu_1,\mu_2)\beta_3,\label{delta}
\end{equation}
where
\begin{eqnarray}
\delta_0&=&-\frac{1}{4} \mu_1 \Big(\hat t_{12} \left(3 t_{11}^2 t_{21}+2 t_{11} t_{12} t_{22}+t_{12}^2 t_{21}\right)+\hat t_{14} \left(3 t_{11}^2 t_{21}+2 t_{11} t_{12} t_{22}+t_{12}^2 t_{21}\right)\nonumber\\
&&+(\hat t_{22}+\hat t_{24}) \left(t_{11}^2 t_{22}+2 t_{11} t_{12} t_{21}+3 t_{12}^2 t_{22}\right)\Big)\label{delta_0}\\
\delta_{\alpha}&=&-\frac{3}{8} \left(t_{11}^2+t_{12}^2\right) (\hat t_{12} t_{11}+\hat t_{14} t_{11}+t_{12} (\hat t_{22}+\hat t_{24}))\\
\delta_{\beta}&=&-\frac{3}{8} \left(t_{31}^2+t_{32}^2\right) (\varepsilon (\hat t_{12} t_{31}+\hat t_{14} t_{31}+t_{32} (\hat t_{22}+\hat t_{24}))+\hat t_{14} t_{31}+\hat t_{24} t_{32})\label{delta_beta3}
\end{eqnarray}
where $t_{ij}$ and $\hat t_{ij}$ are the elements of $\mathbf T$ and of its inverse, respectively.

\begin{figure}[t]
\begin{centering}
\includegraphics[trim = 10mm 0mm 0mm 0mm,width=0.32\textwidth]{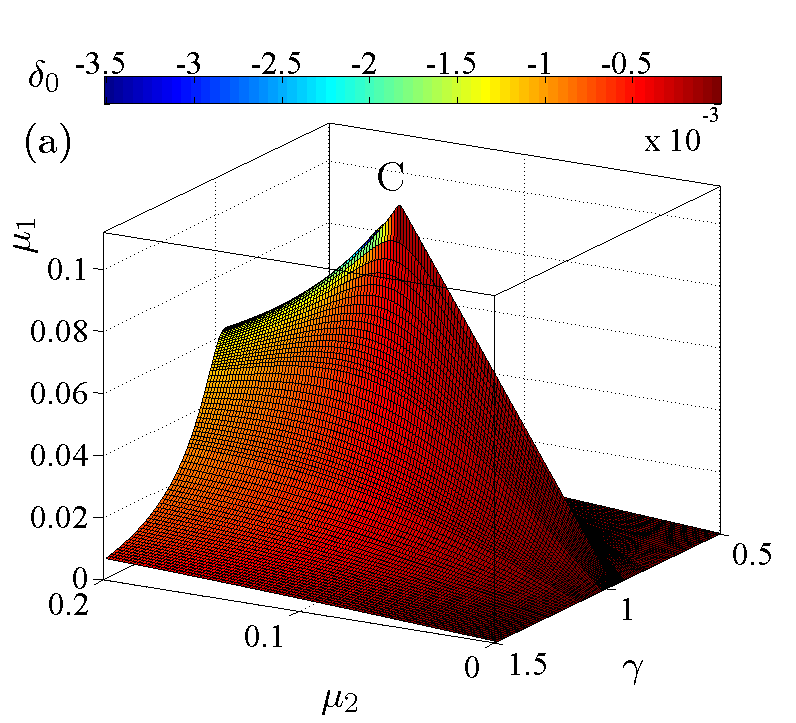}
\includegraphics[trim = 10mm 0mm 0mm 0mm,width=0.32\textwidth]{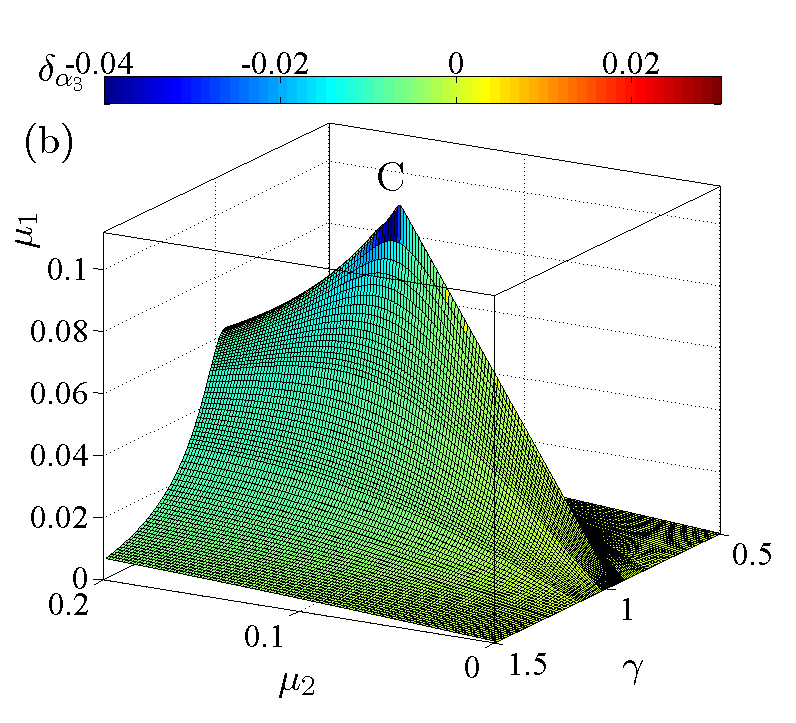}
\includegraphics[trim = 10mm 0mm 0mm 0mm,width=0.32\textwidth]{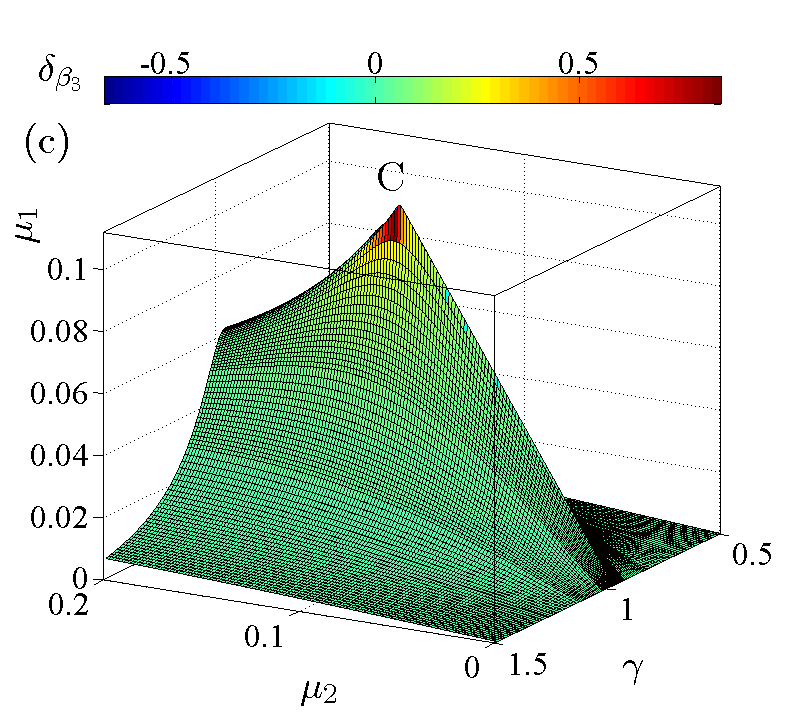}\\
\includegraphics[trim = 10mm 0mm 0mm 0mm,width=0.32\textwidth]{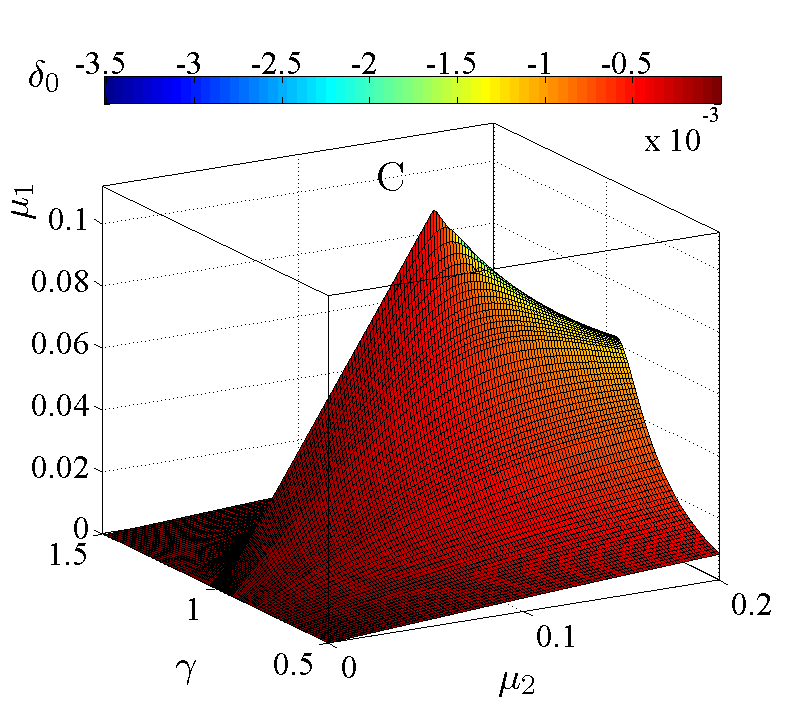}
\includegraphics[trim = 10mm 0mm 0mm 0mm,width=0.32\textwidth]{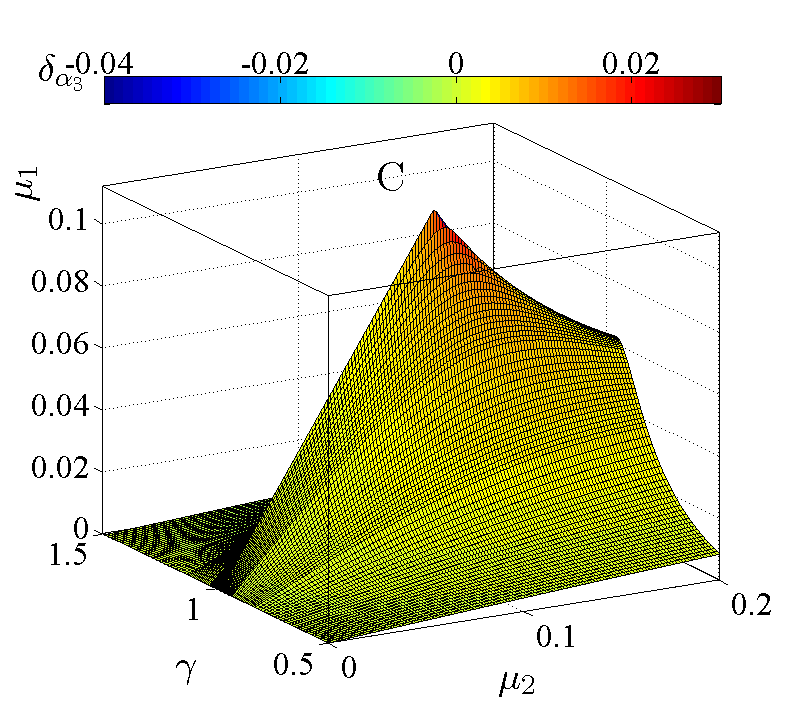}
\includegraphics[trim = 10mm 0mm 0mm 0mm,width=0.32\textwidth]{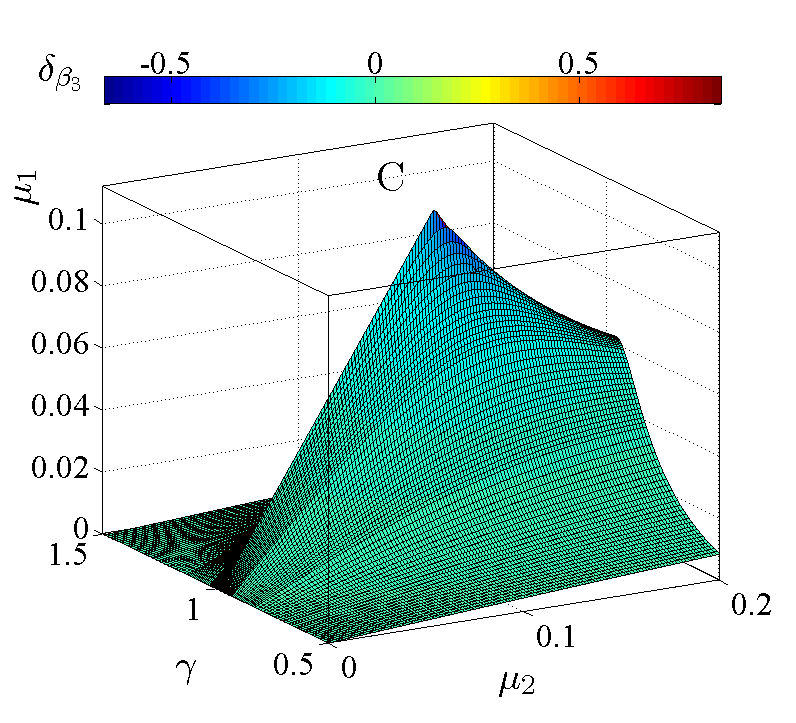}
\par\end{centering}
\caption{\label{coef_delta} Values of (a) $\delta_0$, (b) $\delta_{\alpha}$ and (c) $\delta_{\beta}$ along the stability boundary in the $\mu_1,\mu_2,\gamma$ space for $\varepsilon=0.05$. The color indicates the value of the corresponding coefficient. To facilitate the visualization, two different views of the same surface are given (top and bottom plots).}
\end{figure}
\begin{figure}[!]
\begin{centering}
\includegraphics[trim = 10mm 10mm 10mm 10mm,width=0.32\textwidth]{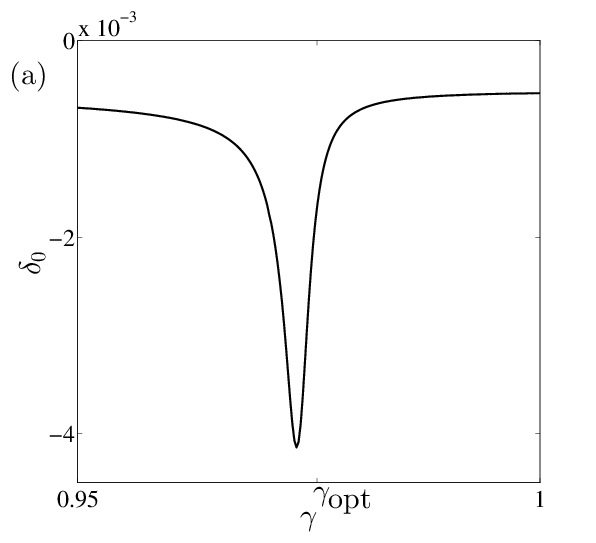}
\includegraphics[trim = 7mm 10mm 13mm 10mm,width=0.32\textwidth]{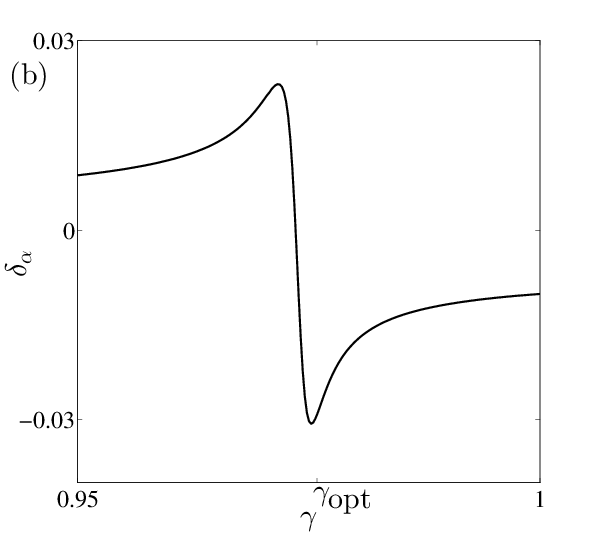}
\includegraphics[trim = 10mm 10mm 10mm 10mm,width=0.32\textwidth]{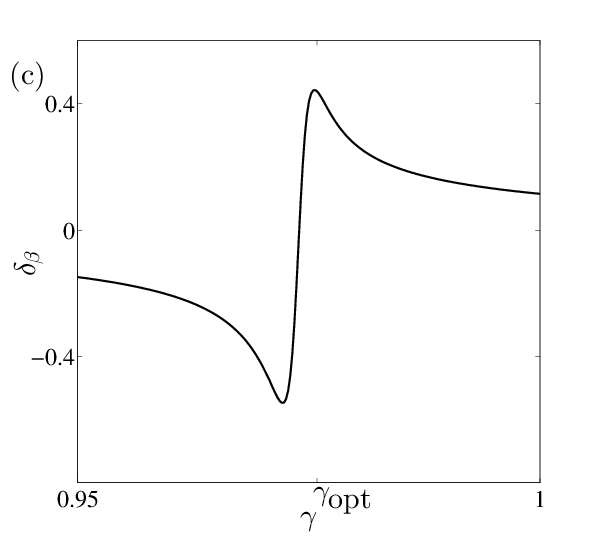}
\par\end{centering}
\caption{\label{coef_delta_mu2fix} Values of (a) $\delta_0$, (b) $\delta_{\alpha}$ and (c) $\delta_{\beta}$ for $\mu_2=0.12$ and $\varepsilon=0.05$.}
\end{figure}

From Eq. (\ref{singH_nf}), we see that the bifurcation is supercritical if $\delta<0$ and subcritical if $\delta>0$. Our objective should therefore be to design the nonlinear spring $\beta_3$ of the NLTVA to impose negative values of $\delta$. 
Ideally, this nonlinear tuning should be carried out in the vicinity of $\gamma_{\text{opt}}$ and $\mu_{2\text{opt}}$ so as to maintain LCO onset at large values of $\mu_1$ (see Section \ref{Elimination}). Fig. \ref{coef_delta} displays the values of the coefficients $\delta_0$, $\delta_{\alpha}$ and $\delta_{\beta}$ along the stability boundary in the $\mu_1,\mu_2,\gamma$ space. Fig. \ref{coef_delta_mu2fix} considers a section of these plots for $\mu_2=0.12$ ($\mu_2=1.1\mu_{2\text{opt}}$) for which stability is lost through a single Hopf bifurcation.

For a system with no structural nonlinearities ($\alpha_3=\beta_3=0$), i.e., a LTVA connected to a classical VdP oscillator, $\delta_0$ is negative, and the system is always supercritical. We now consider a LTVA connected to a VdPD oscillator ($\beta_3=0$), Figs. \ref{coef_delta}(b) and \ref{coef_delta_mu2fix}(b) evidence a symmetric behavior for $\delta_{\alpha}$, i.e., it is positive (negative) below (above) $\gamma \approx \gamma_{\text{opt}}$. This uncertainty on the sign of $\delta_{\alpha}$ in the region of optimum tuning poses an important practical difficulty, because a supercritical bifurcation cannot confidently be enforced in this region. For positive values of $\alpha_3$, the solution to avoid a catastrophic bifurcation scenario is to detune the LTVA toward greater values of $\gamma$, which guarantees negative values of $\delta_{\alpha}$. However, Fig. \ref{mu1_lambda}(c) indicates that this detuning is associated with a significant decrease in the value of $\mu_{1\text{max}}$. For instance, considering $\gamma=1$ decreases $\mu_{1\text{max}}$ by approximately 30\%. Following a similar reasoning, we conclude that the LTVA should be detuned toward smaller values of $\gamma$ for negative values of $\alpha_3$.

The NLTVA presents increased flexibility with respect to the LTVA, because $\beta_3$ represents an additional tuning parameter.
However, Figs. \ref{coef_delta}(c) and \ref{coef_delta_mu2fix}(c) show that the sign of $\delta_{\beta}$ in the optimal tuning region is as difficult to predict as for $\delta_{\alpha}$.
Interestingly, $\delta_{\alpha}$ and $\delta_{\beta}$ have consistently an opposite sign.
It can be verified that $\delta_{\alpha}/\delta_{\beta}\approx-0.05$ close to $\gamma_{\text{opt}}$, thus $\delta_{\alpha}\alpha_3+\delta_{\beta}\beta_3 \approx 0$ if $\beta_3\approx 0.05 \alpha_3$.
In other words, {\it the potentially detrimental effect of the structural nonlinearity of the VdPD on the bifurcation behavior can be compensated through a proper design of the NLTVA's nonlinearity}. Unlike the LTVA, the NLTVA can therefore be designed to enforce supercritical bifurcations in the optimal tuning region.

\subsection{Two intersecting single Hopf bifurcations}

Although the investigation of the double Hopf bifurcation that occurs along the black line in Fig. \ref{stab_3d} is beyond the scope of this paper, the separate analysis of the two intersecting single Hopf bifurcations gives already some insight into the dynamics. The eigenvalues of $\mathbf{W}$ at point C are
$\lambda_{1,2}=\pm j$ and $\lambda_{3,4}=\pm j/\sqrt{1+\varepsilon}$. By performing the analysis outlined in the previous section, first considering $\lambda_{1,2}$ as the critical eigenvalues and then $\lambda_{3,4}$, we obtain respectively
\begin{eqnarray}
\delta_{1,2}&=&\frac{1}{8}\left(-\frac{\varepsilon\sqrt{\varepsilon}}{1+\varepsilon}+\frac{3\sqrt{\varepsilon}}{1+\varepsilon}\alpha_3-\frac{3(1+\varepsilon)}{\sqrt{\varepsilon}}\beta_3\right)\label{delta_1}\\
\delta_{3,4}&=&\frac{3}{8}\left(-\sqrt{\varepsilon}\alpha_3+\frac{(1+\varepsilon)^2}{\sqrt{\varepsilon}}\beta_3\right)\label{delta_2}.
\end{eqnarray}
The remarkable feature of Eqs. (\ref{delta_1}) and (\ref{delta_2}) is that the ratio between $\delta_{\alpha}$ and $\delta_{\beta}$ is constant \begin{equation}
\frac{\delta_{\alpha}}{\delta_{\beta}}=-\frac{\varepsilon}{(1+\varepsilon)^2}. \label{alphapbeta}
\end{equation}
The important practical consequence of this result is that {\it the effect of $\alpha_3$ can be locally entirely compensated by $\beta_3=\varepsilon/(1+\varepsilon)^2\alpha_3$}. It can be shown in fact that Eq. (\ref{alphapbeta}) is valid along the whole line of double Hopf bifurcations. We also note that Eq. (\ref{alphapbeta}) is in excellent agreement with the expression $\delta_{\alpha}/\delta_{\beta}\approx-0.05$ obtained in the single Hopf case for $\varepsilon=0.05$.

\begin{figure}[!]
\begin{centering}
\includegraphics[trim = 10mm 10mm 0mm 10mm,width=0.36\textwidth]{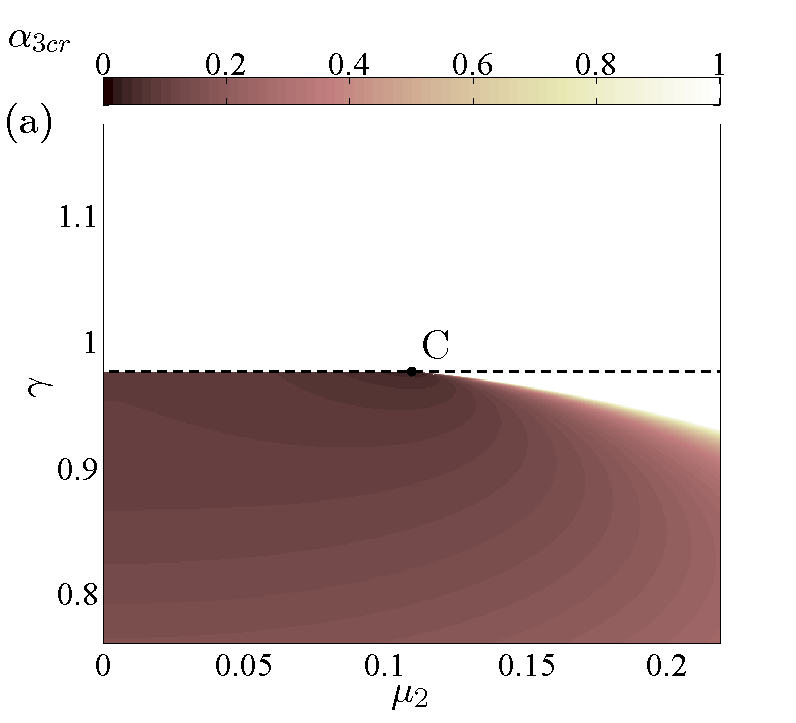}
\includegraphics[trim = 10mm 10mm 0mm 10mm,width=0.36\textwidth]{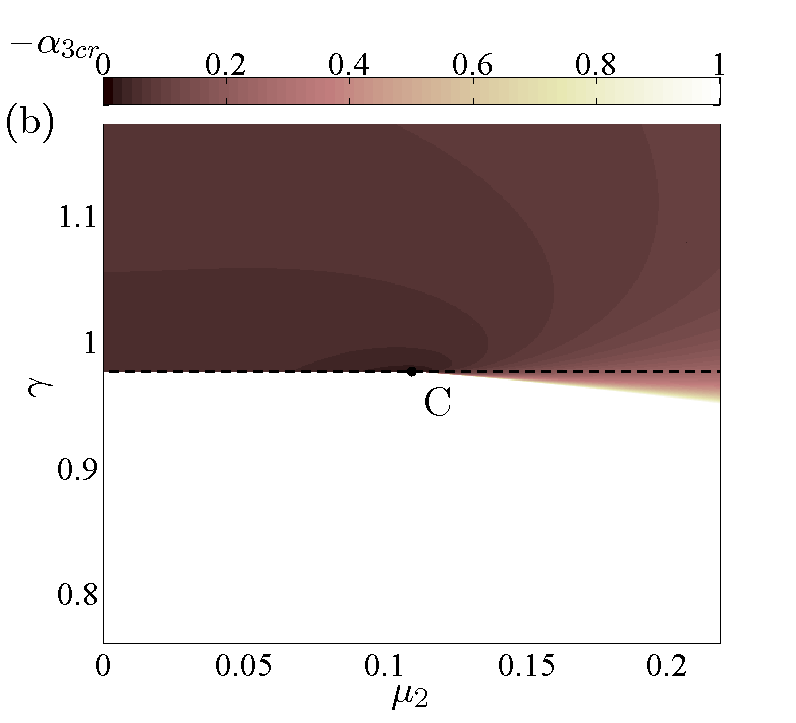}
\par\end{centering}
\caption{\label{max_alpha}Maximal value of $\alpha_3$ below which supercritical bifurcations for the LTVA are guaranteed ($\varepsilon=0.05$, $\beta_3=0$). (a) Positive values of $\alpha_3$; (b) negative values of $\alpha_3$. The dashed line corresponds to $\gamma=\gamma_{\text{opt}}$ and the dot to $\mu_2=\mu_{2\text{opt}}$. }
\end{figure}

\begin{figure}[!]
\begin{centering}
\includegraphics[trim = 0mm 10mm 10mm 10mm,width=0.36\textwidth]{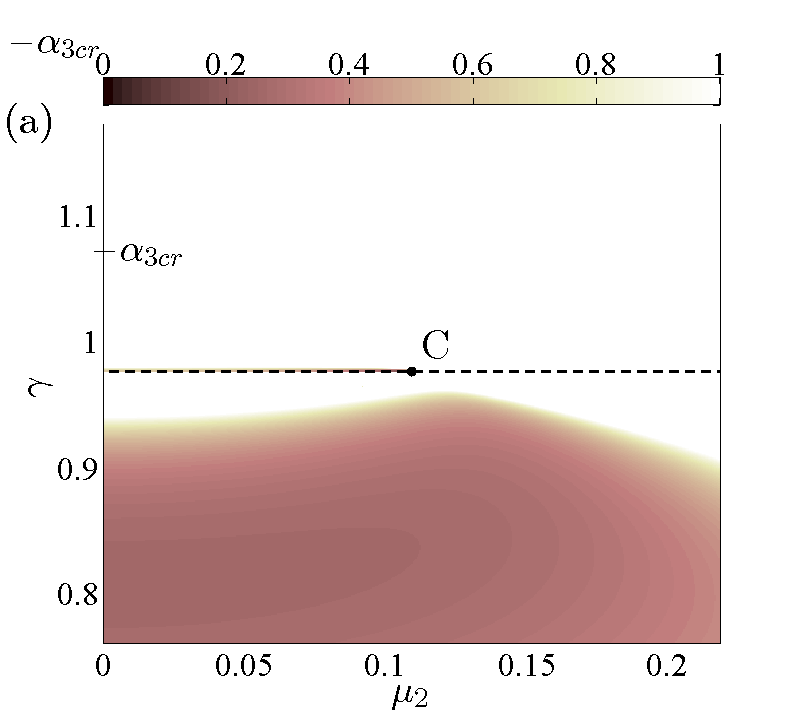}
\includegraphics[trim = 0mm 10mm 10mm 10mm,width=0.36\textwidth]{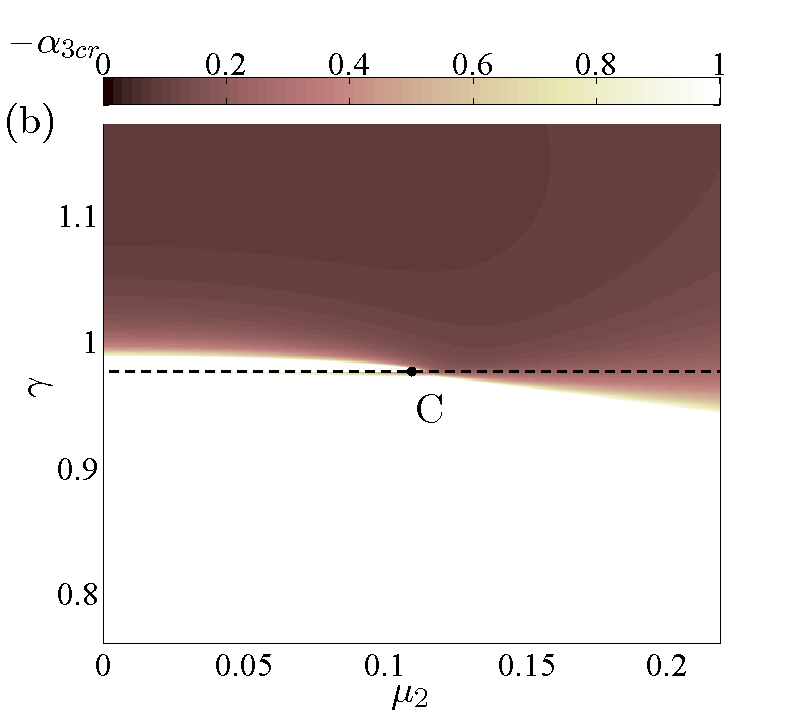}
\par\end{centering}
\caption{\label{max_alpha2}Maximal value of $\alpha_3$ below which supercritical bifurcations for the NLTVA are guaranteed ($\varepsilon=0.05$, $\beta_3=\alpha_3\varepsilon/(1+\varepsilon)^2$). (a) Positive values of $\alpha_3$; (b) negative values of $\alpha_3$. The dashed line corresponds to $\gamma=\gamma_{\text{opt}}$ and the dot to $\mu_2=\mu_{2\text{opt}}$.}
\end{figure}

Although Eq. (\ref{alphapbeta}) is exactly valid for the double Hopf bifurcations, it is only approximate in other cases.
Thus, choosing $\beta_3=\varepsilon/(1+\varepsilon)^2\alpha_3$ does not necessarily guarantee that the system will lose stability through a supercritical Hopf bifurcation.
Considering that $\delta=\delta_0+\delta_{\alpha}\alpha_3+\delta_{\beta}\beta_3$ and imposing $\beta_3=\varepsilon/(1+\varepsilon)^2\alpha_3$, it is possible to calculate the value of $\alpha_3$ such that $\delta$ becomes positive, i.e., the critical value of $\alpha_3$ for which the absorber fails to enforce supercriticality.
This value is given by $\alpha_{3\text{cr}}=-\delta_0/\delta_{\alpha}$ in the case of a LTVA ($\beta_3=0$) and by $\alpha_{3\text{cr}}=-\delta_0/(\delta_{\alpha}+\delta_{\beta}\varepsilon/(1+\varepsilon)^2)$ for a NLTVA with the proposed tuning rule (\ref{alphapbeta}).
The critical values of $\alpha_3$ calculated in both cases are illustrated in Figs. \ref{max_alpha} and \ref{max_alpha2}, respectively.
To avoid very large values in the vicinity of double Hopf bifurcations (where $\alpha_{\text{cr}}\to\infty$), the color maps were trimmed at 1.

Considering the LTVA (Fig. \ref{max_alpha}), it is seen that the point C of optimal tuning of $\gamma$ and $\mu_2$ lies at the boundary between 0 and 1, resulting in a design with virtually zero robustness. For positive (negative) values of $\alpha_3$, the solution for a robust absorber is to increase either $\gamma$ or $\mu_2$ (decrease $\gamma$), which necessarily results in an earlier LCO onset, i.e., $\mu_{1\text{max}}<\sqrt{\varepsilon}/2$.
If a NLTVA is adopted (Fig. \ref{max_alpha2}), for positive values of $\alpha_3$, point C lies well inside the region where supercriticality is guaranteed, which clearly highlights the benefit of the NLTVA. For negative values of $\alpha_3$, the optimal point lies close to the boundary between 0 and 1, which means that there is much less margin for a robust design than for positive $\alpha_3$.
However, compared to the LTVA with negative $\alpha_3$, the NLTVA still possesses a larger region of supercritical behavior. Specifically, there is a new region $\gamma\approx\gamma_{\text{opt}}$ and $\mu_2<\mu_{2\text{opt}}$ in which supercriticality can be guaranteed.

\begin{figure}[!]
\begin{centering}
\includegraphics[trim = 10mm 10mm 10mm 10mm,width=0.4\textwidth]{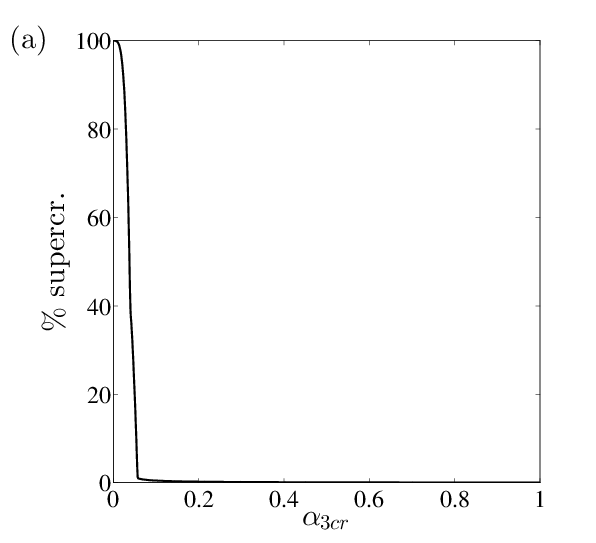}
\includegraphics[trim = 10mm 10mm 10mm 10mm,width=0.4\textwidth]{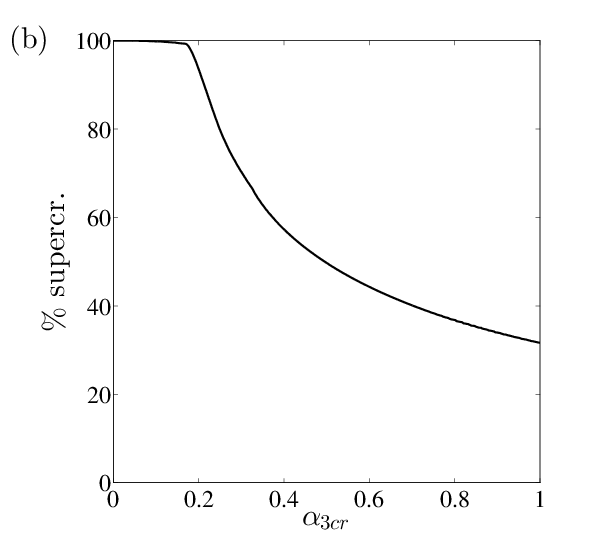}
\par\end{centering}
\caption{\label{percentage}Probability to have a supercritical bifurcation for different values of $\alpha_3$ (in absolute value). (a) LTVA; (b) NLTVA. $\gamma$ and $\mu_2$ are within $\pm$ 1\% and $\pm$ 5\% of the corresponding optimum value, respectively.}
\end{figure}

Finally, Fig. \ref{percentage} presents the probability to have a supercritical bifurcation as a function of $\alpha_3$ (in absolute value). To reflect a realistic design scenario, uncertainty of $\pm$ 1\% and $\pm$ 5\% on the values of $\gamma$ and $\mu_2$ around point C, respectively, are considered. Again, the superiority of the NLTVA over the LTVA is evident in these plots.

\section{Reduction of the Amplitude of Limit Cycle Oscillations}\label{Reduction}

At this stage, the linear and nonlinear parameters of the NLTVA have been designed through stability (Section \ref{Elimination}) and bifurcation (Section \ref{BifurcationAnalysis}) analyses, respectively. There is therefore no much freedom left to mitigate the amplitudes of the LCOs in the post-bifurcation regime.
Nevertheless, their amplitude in the vicinity of the stability boundary can be investigated through the adopted local analysis.

Considering the normal form of a Hopf bifurcation, the amplitude of the generated LCOs is proportional to $\sqrt{-k_1/\delta}$, where $k_1$ is the real part of the eigenvalue related to the bifurcation (see Section \ref{SHBIF}). Because $k_1=0$ at the loss of stability, we consider its linear approximation, i.e., $k_1\approx \left(\text dk_1/\text d\mu_1\right)|_{\mu_1=\mu_{1\text{cr}}}(\mu_1-\mu_{1\text{cr}})$. The LCO amplitude in the vicinity of the loss of stability is therefore \begin{equation}
r\approx\sqrt{-\frac{\text dk_1}{\text d\mu_1}\bigg|_{\mu_1=\mu_{1\text{cr}}}\frac{\mu_1-\mu_{1\text{cr}}}{\delta}},\label{r_ampl}
\end{equation}

\begin{figure}
\begin{centering}
\includegraphics[trim = 5mm 10mm 5mm 10mm,width=0.5\textwidth]{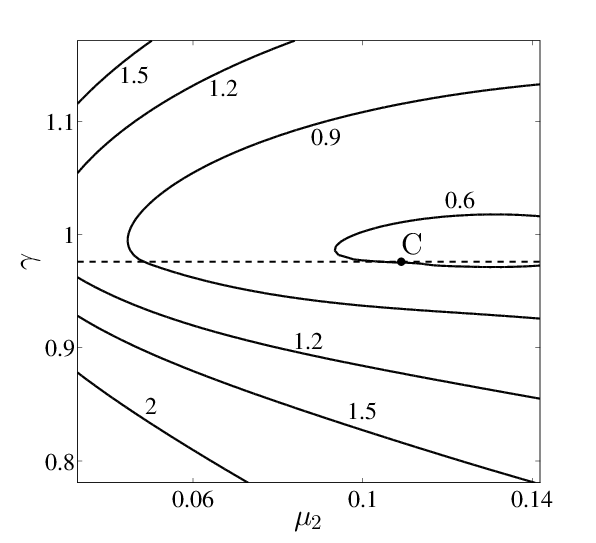}
\par\end{centering}
\caption{\label{q1}LCO iso-amplitude curves along the stability boundary for $\alpha_3=0.08$, $\varepsilon=0.05$, $\beta_3=\varepsilon/(1+\varepsilon)^2\alpha_3=0.0036$ and $\mu_1-\mu_{1\text{cr}}=0.01$.}
\end{figure}

The maximal value of the LCO in physical space is computed by considering $\varphi=\tan^{-1}\left(t_{12}/t_{11}\right)$ in
$q_1\approx t_{11}r\cos\left({\varphi}\right)+t_{12}r\sin\left({\varphi}\right)$,
where $t_{11}$ and $t_{12}$ are the coefficients of $\mathbf T$ as expressed in Eq. (\ref{T}) ($t_{13}$ and $t_{14}$ do not need to be taken into account since $y_3\approx y_4\approx0$).
Iso-amplitude curves of the LCOs along the stability boundary are represented in Fig. \ref{q1}.
It can be observed that the amplitude of the LCOs is minimized close to the optimal tuning region, which signifies that the design of the previous sections is also relevant for LCO mitigation.

\section{Numerical validation of the analytical developments}\label{Validation}

\subsection{Proposed tuning rules for the NLTVA}

In the previous sections, a procedure to optimize the parameters of the NLTVA was proposed, namely \begin{equation}
\gamma_{\text{opt}}=\frac{1}{\sqrt{1+\varepsilon}},\quad \mu_{2\text{opt}}=\frac{1}{2}\sqrt{\frac{\varepsilon}{1+\varepsilon}},\quad\beta_3=\frac{\varepsilon}{(1+\varepsilon)^2}\alpha_3.
\end{equation}
We recall that the 'principle of similarity' was adopted for selecting the mathematical form of the NLTVA nonlinearity.
For these parameters, the system will lose stability at $\mu_1=\sqrt{\varepsilon}/4$ through a supercritical Hopf bifurcation. Because the double Hopf scenario was not fully analyzed and because this scenario would probably result in more involved dynamics, it is probably preferable in practice to detune $\mu_2$ toward values slightly greater than $\mu_{2\text{opt}}$, in order to ensure a supercritical single Hopf bifurcation. Furthermore, Fig. \ref{q1} shows that the LCOs have smaller amplitude for $\mu_2>\mu_{2\text{cr}}$ rather than for $\mu_2<\mu_{2\text{cr}}$.

\subsection{Local and global bifurcation analysis}

\begin{figure}[!]
\begin{centering}
\includegraphics[trim = 10mm 10mm 10mm 10mm,width=0.4\textwidth]{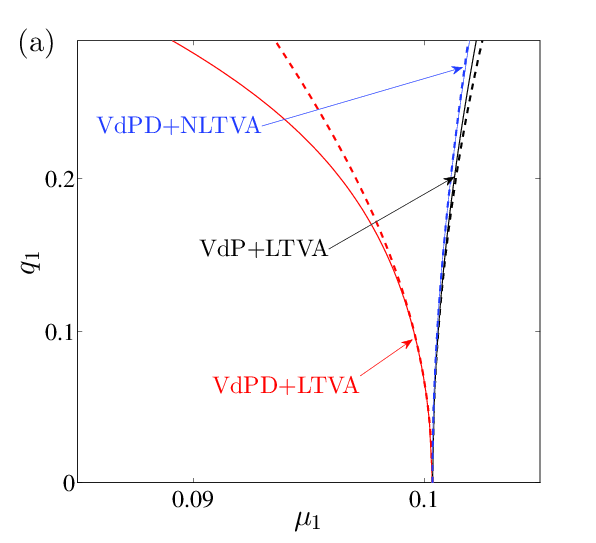}
\includegraphics[trim = 10mm 10mm 10mm 10mm,width=0.4\textwidth]{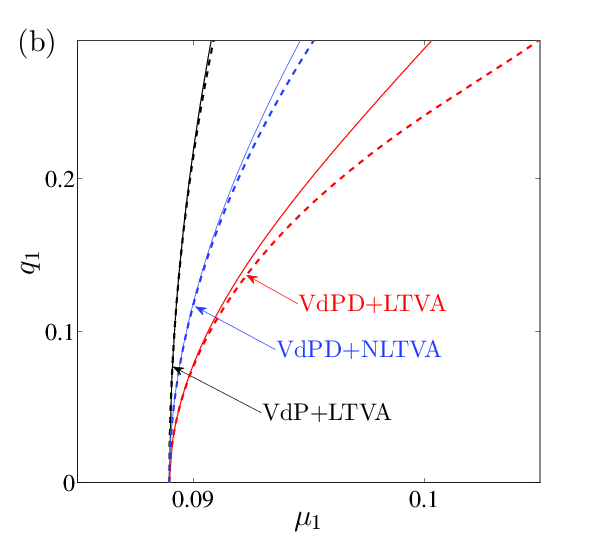}
\par\end{centering}
\caption{\label{bif_singlHopf}Bifurcation diagrams for $\mu_2=0.12$ and $\varepsilon=0.05$. (a) $\gamma=0.970$, (b) $\gamma=0.985$. VdP+LTVA: $\alpha_3=0$ and $\beta_3=0$; VdPD+LTVA: $\alpha_3=0.3$ and $\beta_3=0$; VdPD+NLTVA: $\alpha_3=0.3$ and $\beta_3=0.0136$. Thick dashed lines: numerical results (MATCONT); thin solid lines: analytical results.}
\end{figure}

Bifurcation diagrams predicted using the analytical developments of Section \ref{SHBIF} and the numerical continuation software MATCONT \cite{Dhooge} are depicted in Fig. \ref{bif_singlHopf}. Slightly detuned linear parameters, i.e., $\mu_2=0.12$, $\gamma=0.97/0.985$, are considered to show the robustness of our findings. Loss of stability occurs through a single Hopf bifurcation for the two parameter sets. Fig. \ref{bif_singlHopf} presents an excellent qualitative agreement between the analytical and numerical curves; the quantitative differences observed at higher values of $q_1$ are due to the fact that the analytical results are only valid locally. When there is no structural nonlinearity ($\alpha_3=0$ and $\beta_3=0$), the bifurcation remains supercritical, and the LTVA works effectively on the classical VdP oscillator. The introduction of the structural nonlinearity ($\alpha_3=0.3$) in the VdP oscillator gives rise to a subcritical or supercritical bifurcation in Figs. \ref{bif_singlHopf}(a) and (b), respectively. This result confirms the difficulty to predict the bifurcation behavior of the coupled VdPD and LTVA system in the optimal region; it also highlights the detrimental role played by the structural nonlinearity of the VdPD oscillator. Conversely, the introduction of nonlinearity in the absorber ($\beta_3=\varepsilon/(1+\varepsilon)^2\alpha_3=0.0136$) allows to guarantee a supercritical bifurcation, as for the system without nonlinearity. The compensation effect brought by the NLTVA is therefore clearly demonstrated.

\begin{figure}[!]
\begin{centering}
\includegraphics[trim = 15mm 10mm 5mm 10mm,width=0.4\textwidth]{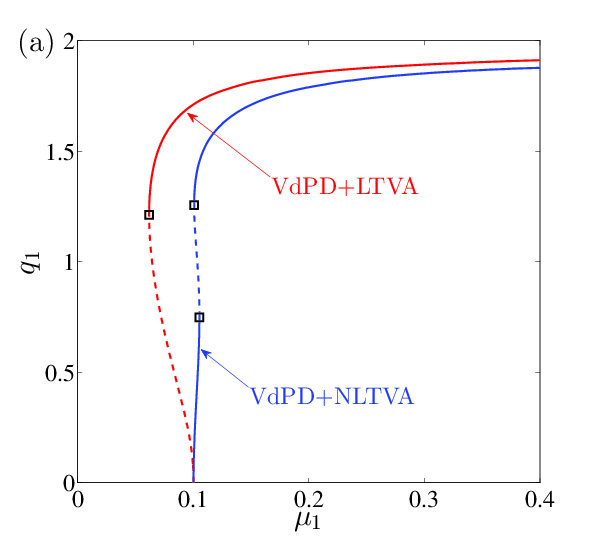}
\includegraphics[trim = 10mm 10mm 10mm 0mm,width=0.4\textwidth]{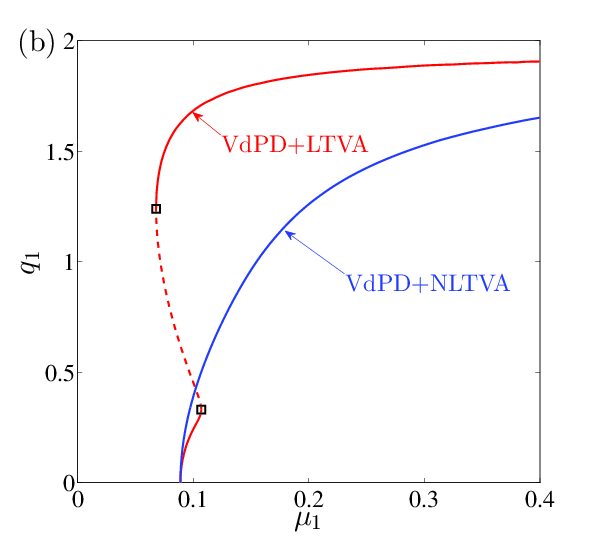}
\par\end{centering}
\caption{\label{bif_singlHopf_large}Bifurcation diagrams for $\mu_2=0.12$, $\alpha_3=0.3$ and $\varepsilon=0.05$. (a) $\gamma=0.970$, (b) $\gamma=0.985$. VdPD+LTVA: $\beta_3=0$; VdPD+NLTVA: $\beta_3=0.0136$. Squares: fold bifurcations. }
\end{figure}

The analytical developments are valid only in the neighborhood of the bifurcation leading to LCO onset. The MATCONT software \cite{Dhooge} is now utilized to investigate large-amplitude LCOs.
Fig. \ref{bif_singlHopf_large} plots bifurcation diagrams for the case of a single Hopf bifurcation at the loss of stability. The same parameter values as those in Fig. \ref{bif_singlHopf} are used, but greater values of $q_1$ are investigated. A major  difference with the local analysis is that fold bifurcations that can turn supercritical behavior into subcritical behavior are now encountered, which confirms the importance of global analysis. If the pair of folds that appears for the NLTVA in Fig. \ref{bif_singlHopf_large}(a) cannot be considered as particularly detrimental, this is not the case for the LTVA in Fig. \ref{bif_singlHopf_large}(b), where bistability in a significant portion of the stable region compromises the robustness of the linear absorber.

\begin{figure}[!]
\begin{centering}
\includegraphics[trim = 15mm 10mm 5mm 10mm,width=0.4\textwidth]{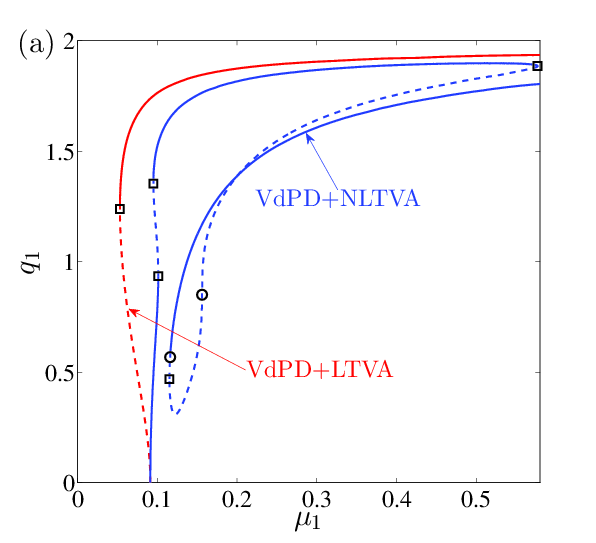}
\includegraphics[trim = 10mm 10mm 10mm 0mm,width=0.4\textwidth]{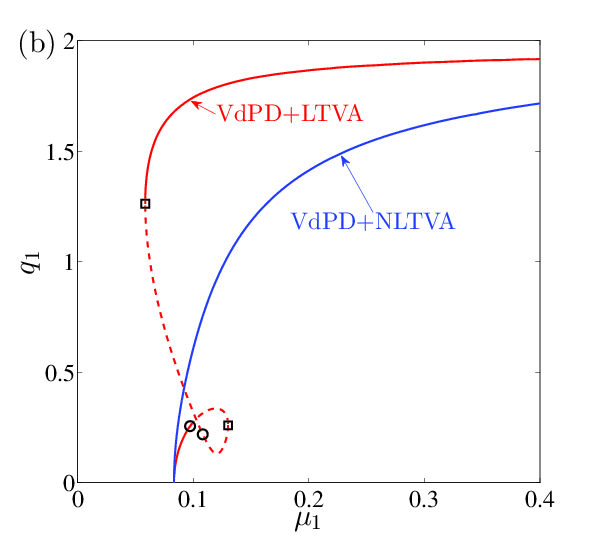}
\par\end{centering}
\caption{\label{bif_doubleHopf_large}Bifurcation diagrams for $\mu_2=0.097$, $\alpha_3=0.3$ and $\varepsilon=0.05$. (a) $\gamma=0.970$, (b) $\gamma=0.985$. VdPD+LTVA: $\beta_3=0$; VdPD+NLTVA: $\beta_3=0.0136$. Squares: fold bifurcations; circles: secondary Hopf bifurcations. }
\end{figure}

Fig. \ref{bif_doubleHopf_large} represents the same results for a lower value of $\mu_2$ for which a double Hopf bifurcation is expected. The bifurcation diagrams are more complex with secondary Hopf (or Neimark-Sacker) bifurcations observed both for the LTVA and NLTVA; their analysis is beyond the scope of this paper. Apart from these new bifurcations, we note that the general trend of the curves is similar to that in Fig. \ref{bif_singlHopf_large}, demonstrating a certain robustness of the absorbers with respect to parameter variations.

\begin{figure}[!]
\begin{centering}
\includegraphics[trim = 15mm 10mm 5mm 10mm,width=0.4\textwidth]{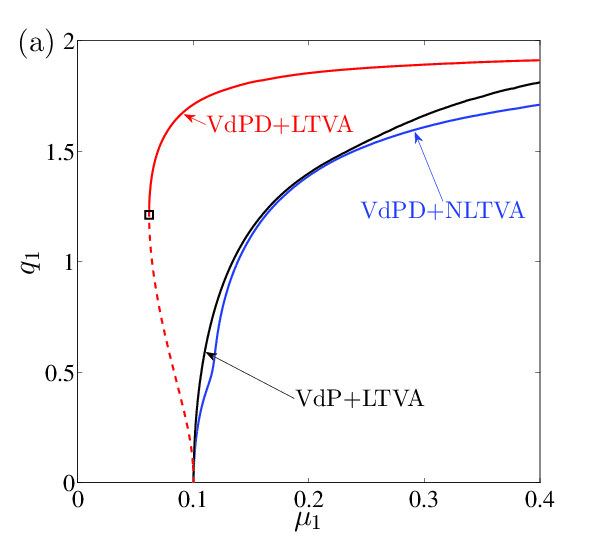}
\includegraphics[trim = 10mm 10mm 10mm 0mm,width=0.4\textwidth]{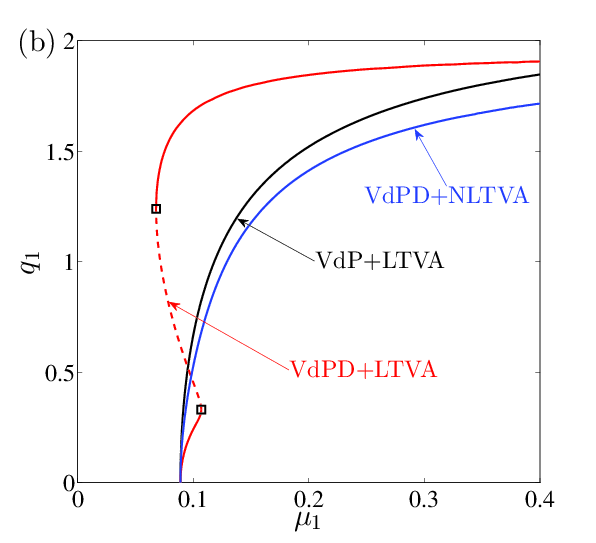}
\par\end{centering}
\caption{\label{bif_singlHopf_large2}Bifurcation diagrams for $\mu_2=0.12$ and $\varepsilon=0.05$. (a)
$\gamma=0.970$, (b) $\gamma=0.985$. VdP+LTVA: $\alpha_3=0$ and $\beta_3=0$; VdPD+LTVA: $\alpha_3=0.3$ and $\beta_3=0$; VdPD+NLTVA: $\alpha_3=0.3$ and $\beta_3=0.018$. Squares: fold bifurcations. }
\end{figure}

Fig. \ref{bif_singlHopf_large2} considers again the single Hopf case, but with a slightly greater value of the nonlinear coefficient of the NLTVA, i.e., $\beta_3=0.018$ instead of $0.0136$. The NLTVA clearly outperforms the LTVA: not only LCO amplitudes are significantly smaller, but there is the complete absence of dangerous bistable regions. Another important result is the strong resemblance between the behaviors of the VdP+LTVA and VdPD+NLTVA systems. This suggests that the compensation of the nonlinearity of the VdPD by the nonlinearity of the NLTVA, which was observed in previous sections, is also valid at larger amplitudes.

\subsection{Validation of the principle of similarity}

\begin{figure}[!]
\begin{centering}
\includegraphics[trim = 15mm 10mm 5mm 10mm,width=0.4\textwidth]{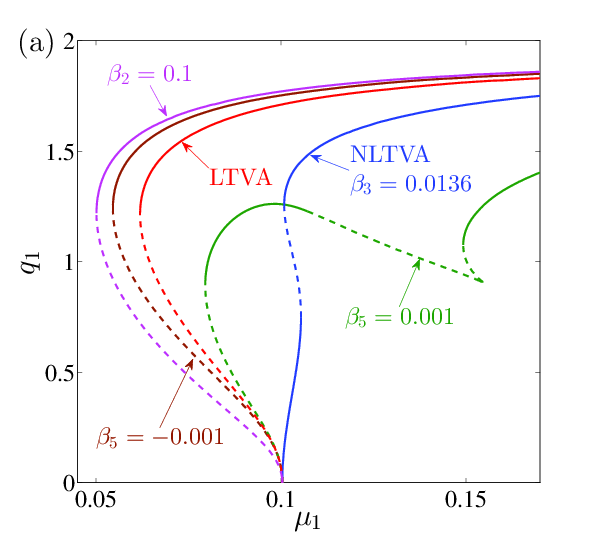}
\includegraphics[trim = 10mm 10mm 10mm 0mm,width=0.4\textwidth]{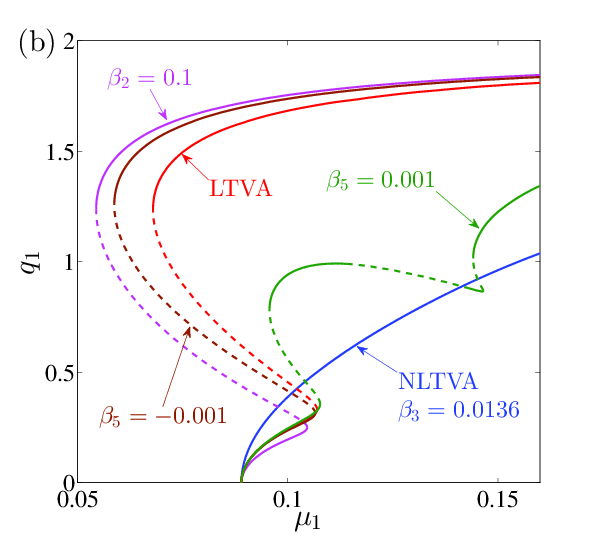}
\par\end{centering}
\caption{\label{235}Bifurcation diagrams for $\mu_2=0.12$, $\varepsilon=0.05$ and $\alpha_3=0.03$ (VdPD primary system) with an attached linear (LTVA), quadratic ($\beta_2$), cubic (NLTVA) or quintic ($\beta_5$) absorber. (a)
$\gamma=0.970$, (b) $\gamma=0.985$.
Solid lines: stable, dashed lines: unstable. }
\end{figure}

In Section \ref{BifurcationAnalysis}, the nonlinear spring of the NLTVA was chosen to be cubic according to the 'principle of similarity'.
Absorbers with quadratic ($\beta_2$)  and quintic ($\beta_5$)  nonlinearities are also considered in this section. For $\gamma=0.97$ in Fig. \ref{235}(a), the NLTVA is the only absorber which loses stability through a supercritical bifurcation. For greater amplitudes, unlike the other absorbers, the NLTVA does not exhibit any bistability in the region where the trivial solution is stable.
However, besides the NLTVA, the absorber with quintic nonlinearity also exhibits a beneficial effect with respect to LCOs amplitude.
For $\gamma=0.985$ in Fig. \ref{235}(b), all absorbers lose stability through a supercritical bifurcation, but the NLTVA is the only absorber that does not possess a fold bifurcation that generates subcriticalities at greater amplitudes. The remaining LCOs for the NLTVA have also a much smaller amplitudes. All the results clearly validate the principle of similarity for the passive control of self-excited oscillations.

\subsection{Comparison with NES}

As stated in the introduction, the NES has been widely utilized in the last decade for the passive absorption of forced and self-excited oscillations.
The main difference between the NES and the NLTVA, is that the former absorber possesses an essentially nonlinear spring, i.e., without a linear component.
Domany and Gendelman \cite{Domany} showed that the stability boundary of a VdPD oscillator with an attached NES can be approximated by the curve $2\mu_1/\varepsilon=\Lambda/(\Lambda^2+1)$, where $\Lambda=c_2/(m_2\omega_{n1})$.
Thus, the maximal value of $\mu_1$ in order to have stability is obtained for $\Lambda\approx1$ and it is $\mu_{1\text{max,NES}}\approx\varepsilon/4$, which is the square of the value obtainable with the NLTVA ($\mu_{1\text{max}}=\sqrt{\varepsilon}/2$).
For $\varepsilon=0.05$, $\mu_{1\text{max}}$ is therefore approximately 10 times greater for the NLTVA than for the NES.
This demonstrates that the presence of a linear spring in the absorber is critical for the enlargement of the stable region of the trivial solution.

\begin{figure}[!]
\begin{centering}
\includegraphics[trim = 10mm 10mm 5mm 10mm,width=0.32\textwidth]{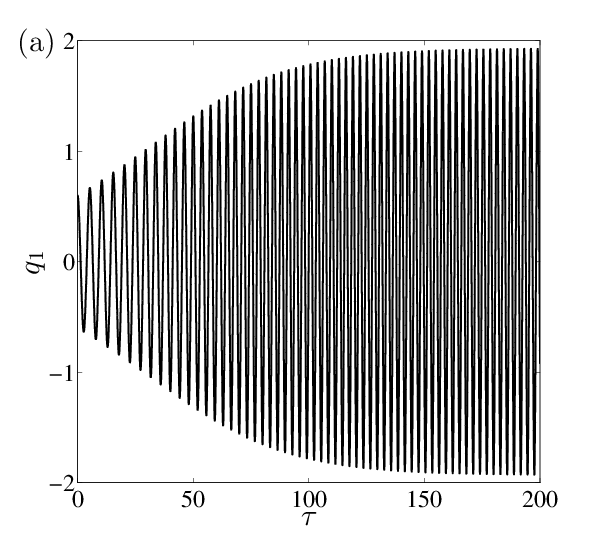}
\includegraphics[trim = 10mm 10mm 5mm 10mm,width=0.32\textwidth]{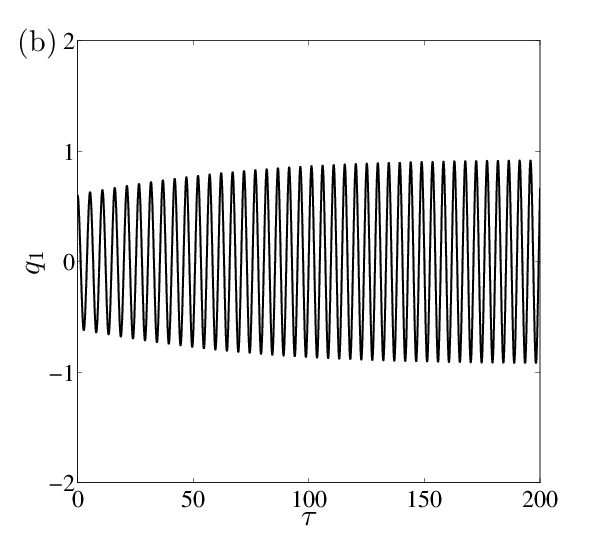}
\includegraphics[trim = 10mm 10mm 5mm 10mm,width=0.32\textwidth]{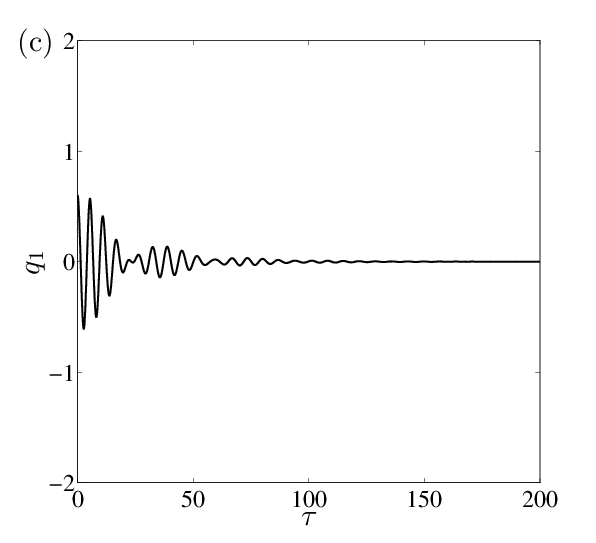}
\par\end{centering}
\caption{\label{fig7}Time series for $\mu_1=0.025$ and $\alpha_3=4/3$. (a) No absorber; (b) NES with $\Lambda=1$ and a cubic spring $\beta_3=0.5333$; (c) NLTVA with optimal linear parameters and a cubic spring $\beta_3=0.0605$.}
\end{figure}

For illustrating the respective performance, a VdPD oscillator with $\mu_1=0.025$ and $\alpha_3=4/3$ is considered. Fig. \ref{fig7}(a) depicts its LCO when no absorber is attached. Figs. \ref{fig7}(b,c) show that the NES can only reduce the LCO amplitude for these parameters whereas the NLTVA can completely eliminate the limit cycle and make the coupled system converge toward the trivial solution.

\begin{figure}[!]
\begin{centering}
\includegraphics[trim = 15mm 10mm 5mm 10mm,width=0.4\textwidth]{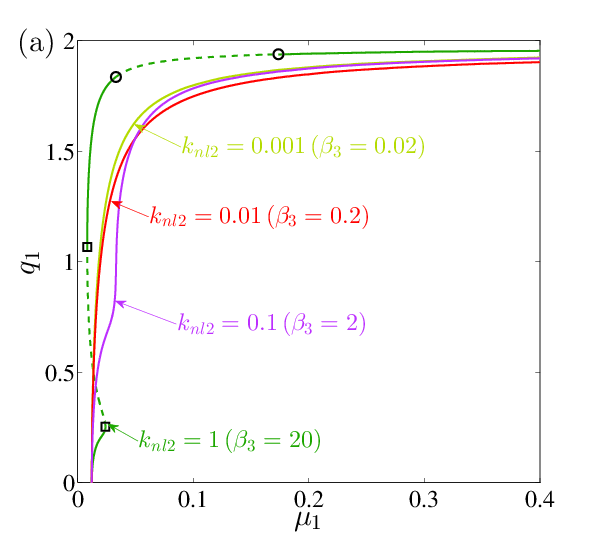}
\includegraphics[trim = 10mm 10mm 10mm 0mm,width=0.4\textwidth]{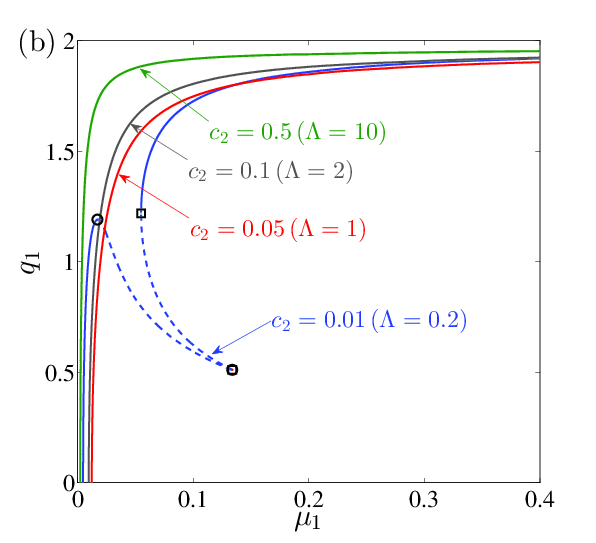}
\par\end{centering}
\caption{\label{bif_NES}Bifurcation diagrams of the VdPD oscillator with an attached NES for $\alpha_3=0.3$ and $\varepsilon=0.05$. (a) $c_2=0.05$ ($\Lambda=1$); (b) $k_{nl2}=0.01$. Squares: fold bifurcations; circles: secondary Hopf bifurcations. }
\end{figure}

The performance of an NES attached to the VdPD oscillator was also investigated using the MATCONT sofware. In Fig. \ref{bif_NES}(a), the NES damper coefficient $c_2$ is optimized with respect to stability and $k_{nl2}$ varies, whereas $k_{nl2}=0.01$ and $c_2$ varies in Fig. \ref{bif_NES}(b). Although a detailed performance analysis is beyond the scope of this paper, variation of $k_{nl2}$ and $c_2$ are unable to improve significantly the NES performance. The comparison between Figs. \ref{bif_singlHopf_large2} and \ref{bif_NES} evidences that, besides a much larger stable region, the NLTVA generates smaller LCO amplitudes than the NES.

\section{Conclusions}\label{Conclusions}

The purpose of this paper was to investigate the performance of the NLTVA for suppression of self-excited oscillations of mechanical systems. A distinct feature of this absorber is that the mathematical form of its nonlinearity is selected according to a principle of similarity with the nonlinear primary system. Thanks to detailed stability and bifurcation analyses, a complete analytical design of the NLTVA was obtained, which was validated by detailed numerical calculations in the MATCONT software. 

Thanks to the complementary roles played by the linear and nonlinear springs of the NLTVA, this absorber was shown to be effective for LCO suppression and mitigation, as it maximizes the stability of the trivial equilibrium point, guarantees supercritical bifurcations and reduces the amplitude of the remaining LCOs. 

\vspace{1cm}
\textbf{Competiting interests}. The authors declare no conflicts of interest.

\textbf{Data accessibility}. There is no data related to this research.

\textbf{Authors' contributions}. GH performed the stability and bifurcation analyses, and defined the tuning rule for nonlinearity compensation. GK proposed the idea of the nonlinear tuned vibration absorber, supervised the research, and established the logical organization of the paper.


\textbf{Funding statement}. The work was supported by the European Union (ERC Starting Grant NoVib 307265).

\end{document}